\documentclass[11pt]{amsart}
\usepackage[margin=1in]{geometry}
\usepackage{amsmath,amssymb,amsthm}
\usepackage{hyperref}
\hypersetup{colorlinks=true,linkcolor=blue,citecolor=blue,
pdftitle={Weighted Euler transforms, real-rootedness, and gamma-positivity},
pdfauthor={Yinuo Cheng and Guoce Xin},
pdfkeywords={Weighted Euler transforms, real-rootedness, gamma-positivity,
multiplier sequences, generalized Laguerre polynomials, Ehrhart polynomials}}

\theoremstyle{plain}
\newtheorem{theorem}{Theorem}[section]
\newtheorem{lemma}[theorem]{Lemma}
\newtheorem{proposition}[theorem]{Proposition}
\newtheorem{corollary}[theorem]{Corollary}
\theoremstyle{definition}
\newtheorem{definition}[theorem]{Definition}
\newtheorem{remark}[theorem]{Remark}

\newcommand{\lc}{\operatorname{lc}}

\title{Weighted Euler transforms, real-rootedness, and $\gamma$-positivity}
\author{Yinuo Cheng}
\address{School of Mathematical Sciences, Capital Normal University,
Beijing 100048, P.R.~China}
\email{2250502145@cnu.edu.cn}

\author{Guoce Xin}
\address{School of Mathematical Sciences, Capital Normal University,
Beijing 100048, P.R.~China}
\email{guoce\_xin@163.com}

\subjclass[2020]{Primary 05A15; Secondary 26C10, 52B20, 33C45}
\keywords{Weighted Euler transforms, real-rootedness, $\gamma$-positivity,
multiplier sequences, generalized Laguerre polynomials, Ehrhart polynomials}
\date{}

\begin{document}

\begin{abstract}
For \(a>0\) and a polynomial \(F\) of degree \(m\), let \(H_{a,F}\)
denote the weighted Euler numerator of \(F\), defined as the numerator
of the generating function \(\sum_{k\ge0}(a)_kF(k)x^k/k!\) with
denominator \((1-x)^{a+m}\). For every \(a>0\), we show the following
preservation result: if the numerators \(H_{a,F}\) and \(H_{a,G}\)
have only real nonpositive zeros, then the same holds for \(H_{a,FG}\). The proof relies on an explicit
multiplication formula in falling-factorial coordinates and a
characterization of multiplier sequences for generalized Laguerre bases
due to Br\"and\'en and Ottergren. Furthermore,
when the symmetric degrees coincide with the ordinary degrees of the
input polynomials, the same operation preserves \(\gamma\)-positivity.
For positive integers \(a=r\), the weight is the Ehrhart polynomial of
the standard simplex \(\Delta_{r-1}\). This yields geometric
preservation theorems for products with a single simplex factor, as
well as applications to Gorenstein polytopes and a Reeve tetrahedron
example.
\end{abstract}

\maketitle
\section{Introduction}
In enumerative combinatorics and Ehrhart theory, one frequently evaluates a polynomial \(F\) at nonnegative integers, or equivalently, computes the generating function of the sequence \(F(k)\). Such a function is rational and its numerator often encodes the core combinatorial data. A natural question is how the numerator behaves under multiplication of the input polynomials.

This paper studies a weighted family of such numerators. Fix \(a>0\), and
let \((a)_k=a(a+1)\cdots(a+k-1)\) denote the rising factorial (with
\((a)_0=1\)). Given a nonzero polynomial \(F\in\mathbb{R}[t]\) of degree
\(m\), define its weighted Euler numerator \(H_{a,F}\) via
\begin{equation}\label{eq:Hdef}
\sum_{k\ge0}\frac{(a)_k}{k!}F(k)x^k
=\frac{H_{a,F}(x)}{(1-x)^{a+m}}.
\end{equation}
By Proposition~\ref{prop:explicit}, \(H_{a,F}\) is a polynomial of degree
at most \(m\). As usual, the power \((1-x)^{-a}\) is understood as a
formal power series, or analytically for \(|x|<1\).

Two special cases explain the role of the parameter \(a\). For \(a=1\)
the weight is simply \(1\), yielding the ordinary generating function of \(F(k)\). For
\(a=r\in\mathbb{N}\), the weight \((r)_k/k!\) is the Ehrhart
polynomial of the standard simplex \(\Delta_{r-1}\); this is the
source of the geometric applications in Section~\ref{sec:geom}.

For the unweighted case \(a=1\), Wagner \cite{Wagner} proved that \(H_{1,FG}\) has only real nonpositive zeros whenever \(H_{1,F}\) and \(H_{1,G}\) do. Br\"and\'en, Ferroni, and Jochemko \cite{BFJ} established further preservation results, including \(\gamma\)-positivity. For general \(a>0\), Karp and Vishnyakova \cite[Open problem 2]{KV} asked whether real negative zeros of \(H_{a,F}\) and \(H_{a,G}\) imply the same for \(H_{a,FG}\). Vishnyakova \cite[Open Problem 1.9]{V} then posed the analogous question for nonpositive zeros. Since \(H_{a,F}(0)=F(0)\), the nonpositive version is more general, as it additionally allows input polynomials with vanishing constant terms. Below, we prove this nonpositive preservation for all \(a>0\) (Theorem~\ref{thm:rr-main}), and deduce the corresponding \(\gamma\)-positivity result (Theorem~\ref{thm:gamma-main}).

For \(m\ge0\), let
\begin{equation}\label{eq:Gammam}
\Gamma_m=\left\{\sum_{i=0}^{\lfloor m/2\rfloor}
\gamma_i x^i(1+x)^{m-2i}:\gamma_i\ge0\right\}.
\end{equation}
A nonzero polynomial \(P\in\Gamma_m\) is said to be \(\gamma\)-positive (or nonnegative in some literature) with
symmetric degree \(m\). Such a polynomial satisfies the symmetry
\(P(x)=x^mP(1/x)\), although its ordinary degree can be strictly
smaller than \(m\). As usual, nonzero constants are regarded as
real-rooted.

We prove two preservation theorems. The first addresses real-rootedness.

\begin{theorem}[Real-rootedness]\label{thm:rr-main}
Let \(a>0\), and let \(F,G\in\mathbb{R}[t]\) be nonzero polynomials.
If every zero of \(H_{a,F}\) and \(H_{a,G}\) is real and nonpositive,
then the same holds for \(H_{a,FG}\).
\end{theorem}

The second addresses \(\gamma\)-positivity, under a subtle degree condition.

\begin{theorem}[\(\gamma\)-positivity]\label{thm:gamma-main}
Let \(a>0\), and let \(F,G\in\mathbb{R}[t]\) be nonzero polynomials
of degrees \(m,n\), respectively. If \(H_{a,F}\in\Gamma_m\) and
\(H_{a,G}\in\Gamma_n\), then \(H_{a,FG}\in\Gamma_{m+n}\).
\end{theorem}

We begin by passing to falling-factorial coordinates. Let \(\mathcal{L}_a(t^{\underline{i}})=(a)_i y^i\) denote the linear change of basis, where \(t^{\underline{i}}\) is the falling factorial. This map converts the multiplication of polynomials into the bilinear operation
\begin{equation}\label{eq:star}
(f\star_ag)(y)=\sum_{k\ge0}\frac{[y(y+1)]^k}{k!(a)_k}
f^{(k)}(y)g^{(k)}(y).
\end{equation}
Under the substitution \(y=x/(1-x)\), the numerator \(H_{a,F}(x)\) transforms into the polynomial \(\mathcal{L}_a(F)(y)\), and its nonpositive zeros map to the zeros of \(\mathcal{L}_a(F)\) in \([-1,0]\). This polynomial \(\mathcal{L}_a(F)(y)\) coincides exactly with the characterizing polynomial in the deep classification of multiplier sequences for generalized Laguerre bases by Br\"and\'en and Ottergren \cite{BO14}. By invoking their characterization as a black box, we bypass several pages of standard technical estimates and obtain Theorem~\ref{thm:rr-main} in a few lines.

For Theorem~\ref{thm:gamma-main}, bilinearity reduces the problem to the \(\gamma\)-basis polynomials \(x^i(1+x)^{m-2i}\). The \(\gamma\)-positivity for the product of any two such basis polynomials is proved in two steps: Theorem~\ref{thm:rr-main} ensures the product has only real nonpositive zeros, and the Br\"and\'en--Gal theorem then yields \(\gamma\)-positivity.

Now we describe two geometric applications of these algebraic results 
when \(a=r\in\mathbb{N}\). Let \(P\) be a \(d\)-dimensional lattice polytope with Ehrhart polynomial \(L_P\), and let \(\Delta_{r-1}\) be the standard simplex. The identity
\begin{equation}\label{eq:ehrhart-ident}
H_{r,L_P}(x)=h^*_{P\times\Delta_{r-1}}(x)
\end{equation}
allows us to restate Theorems~\ref{thm:rr-main} and~\ref{thm:gamma-main} as follows.

\begin{theorem}[Real-rootedness, geometric form]\label{thm:geom-rr}
Let \(P,Q\) be lattice polytopes of dimensions \(d,e\), and let
\(r\ge1\) be an integer. If \(h^*_{P\times\Delta_{r-1}}\) and
\(h^*_{Q\times\Delta_{r-1}}\) are real-rooted, then
\(h^*_{P\times Q\times\Delta_{r-1}}\) is real-rooted.
\end{theorem}

\begin{theorem}[\(\gamma\)-positivity, geometric form]\label{thm:geom-gamma}
Let \(P,Q\) be lattice polytopes of dimensions \(d,e\), and let
\(r\ge1\) be an integer. If
\(h^*_{P\times\Delta_{r-1}}\in\Gamma_d\) and
\(h^*_{Q\times\Delta_{r-1}}\in\Gamma_e\), then
\(h^*_{P\times Q\times\Delta_{r-1}}\in\Gamma_{d+e}\).
\end{theorem}

A key structural feature of these results is that the output retains only a single simplex factor. This connects to the general study of Ehrhart invariants under polytope operations (see Liu \cite{Liu}). Moreover, Theorem~\ref{thm:geom-gamma} naturally applies to Gorenstein polytopes.

The remainder of the paper is organized as follows. Section~\ref{sec:mult} establishes the multiplication formula, Section~\ref{sec:rr} proves the preservation of real-rootedness, Section~\ref{sec:gamma} proves \(\gamma\)-positivity, and Section~\ref{sec:geom} develops the geometric consequences, including a Reeve tetrahedron example.

\section{The multiplication formula for weighted Euler numerators}\label{sec:mult}

This section reduces the product \(H_{a,FG}\) to an explicit bilinear
formula \eqref{eq:starformula}. This formula is the key tool for the preservation arguments in Sections~\ref{sec:rr} and~\ref{sec:gamma}.

\subsection{Falling-factorial coordinates}
To simplify multiplication, we pass to the basis of falling factorials \(t^{\underline{j}}=t(t-1)\dots(t-j+1)\), with \(t^{\underline{0}}=1\). If \(F(t)=\sum_{j=0}^mc_jt^{\underline{j}}\), define
\begin{equation}\label{eq:La}
\mathcal{L}_a(F)(y)=\sum_{j=0}^mc_j(a)_jy^j.
\end{equation}

 Since $(a)_j>0$, \(\mathcal{L}_a\) is a degree-preserving linear isomorphism. 
 For $F\neq0$ of degree $m$,
\begin{equation}\label{eq:lcLa}
\lc(\mathcal{L}_a(F))=(a)_m\lc(F),
\end{equation}
where $\lc$ extracts the  
leading coefficient. 

\begin{proposition}\label{prop:explicit}
For $F=\sum_{j=0}^mc_jt^{\underline{j}}$, we have
\begin{equation}\label{eq:explicit}
\sum_{k\ge0}\frac{(a)_k}{k!}F(k)x^k
=\frac{\sum_{j=0}^mc_j(a)_jx^j(1-x)^{m-j}}{(1-x)^{a+m}}.
\end{equation}
If $m=\deg F$, the numerator and $\mathcal{L}_a(F)$ are related by the mutually inverse substitutions
\begin{equation}\label{eq:Hint}
H_{a,F}(x)=(1-x)^m\mathcal{L}_a(F)\left(\frac{x}{1-x}\right),
\end{equation}
and
\begin{equation}\label{eq:Lint}
\mathcal{L}_a(F)(y)=(1+y)^mH_{a,F}\left(\frac{y}{1+y}\right).
\end{equation}
Furthermore, the endpoint values satisfy
\begin{equation}\label{eq:Hvalues}
H_{a,F}(0)=F(0),\qquad H_{a,F}(1)=(a)_m\lc(F).
\end{equation}
\end{proposition}

\begin{proof}
The binomial series and the identity $(a)_{j+\ell}=(a)_j(a+j)_\ell$ yield
\[
\sum_{k\ge0}\frac{(a)_k}{k!}k^{\underline{j}}x^k
=\sum_{k\ge j}\frac{(a)_k}{(k-j)!}x^k
=(a)_jx^j\sum_{\ell\ge0}\frac{(a+j)_\ell}{\ell!}x^\ell
=\frac{(a)_jx^j}{(1-x)^{a+j}}.
\]
Multiplying by $c_j$ and summing over \(j\) yields \eqref{eq:explicit}.
When \(m=\deg F\), the numerator in \eqref{eq:explicit} is exactly \((1-x)^m\mathcal{L}_a(F)(\frac{x}{1-x})\), which gives \eqref{eq:Hint}.
Substituting \(x=y/(1+y)\) into \eqref{eq:Hint} gives \eqref{eq:Lint}.
Finally, evaluating \eqref{eq:Hint} at \(x=0\) and \(x=1\) yields \(c_0=F(0)\) and \(c_m(a)_m=(a)_m\lc(F)\). This proves \eqref{eq:Hvalues}.
\end{proof}
\subsection{Two auxiliary identities}

\begin{lemma}\label{lem:fallprod}
For $i,j\ge0$,
\begin{equation}\label{eq:fallprod}
t^{\underline{i}}t^{\underline{j}}
=\sum_{\ell=0}^{\min(i,j)}\binom{i}{\ell}\binom{j}{\ell}\ell!\,
t^{\underline{i+j-\ell}}.
\end{equation}
\end{lemma}
\begin{proof}
Since both sides are polynomials in \(t\), it suffices to assume that \(t\) is a sufficiently large integer. For such \(t\), we give a combinatorial proof. The left-hand side counts pairs of ordered lists of lengths \(i\) and \(j\) chosen from \(t\) labels, with no repetition within each list. Classify such a pair by the number \(\ell\) of common labels. There are \(\binom{i}{\ell}\binom{j}{\ell}\ell!\) ways to choose and match these common labels. Once these are fixed, the remaining \(i+j-\ell\) positions must receive distinct labels, yielding \(t^{\underline{i+j-\ell}}\) choices. This completes the proof.
\end{proof}

\begin{lemma}\label{lem:pochsum}
If $I,J$ are nonnegative integers and $b>0$, then
\begin{equation}\label{eq:pochsum}
\sum_{r=0}^{\min(I,J)}\frac{I^{\underline{r}}J^{\underline{r}}}
{r!(b)_r}=\frac{(b+J)_I}{(b)_I}.
\end{equation}
\end{lemma}

\begin{proof}
We use the generalized binomial coefficients \(\binom{u}{v}=u(u-1)\dots(u-v+1)/v!\) for real \(u\) and integers \(v\ge0\). Multiplying the left-hand side by \((b)_I/I!\) gives
\[
\sum_{r=0}^{\min(I,J)}\binom{b+I-1}{I-r}\binom{J}{r}.
\]
This sum is exactly the coefficient of \(z^I\) in the product \((1+z)^{b+I-1}(1+z)^J=(1+z)^{b+I+J-1}\), which evaluates to \(\binom{b+I+J-1}{I}=(b+J)_I/I!\). This is an instance of the Vandermonde convolution. Dividing by \((b)_I/I!\) completes the proof.
\end{proof}

\subsection{The multiplication formula}
For $f,g\in\mathbb{R}[y]$, we define the induced operation
\begin{equation}\label{eq:starDef}
f\star_ag=\mathcal{L}_a(\mathcal{L}_a^{-1}(f)\mathcal{L}_a^{-1}(g)).
\end{equation}
This operation is bilinear, commutative, and associative. The next proposition gives an explicit formula for it.

\begin{proposition}\label{prop:starformula}
For any $a>0$ and $f,g\in\mathbb{R}[y]$, we have
\begin{equation}\label{eq:starformula}
(f\star_ag)(y)=\sum_{k\ge0}\frac{[y(y+1)]^k}{k!(a)_k}
f^{(k)}(y)g^{(k)}(y).
\end{equation}
If $f$ and $g$ are nonzero, the sum is finite with $k\le\min(\deg f,\deg g)$.
\end{proposition}

\begin{proof}
Since both sides are bilinear, it suffices to verify the identity for monomials \(f(y)=y^i\) and \(g(y)=y^j\). Using \(\mathcal{L}_a^{-1}(y^i)=t^{\underline{i}}/(a)_i\) and Lemma~\ref{lem:fallprod}, we compute
\[
y^i\star_ay^j=\sum_{\ell=0}^{\min(i,j)}A_\ell y^{i+j-\ell},\qquad
A_\ell=\binom{i}{\ell}\binom{j}{\ell}\ell!\,
\frac{(a)_{i+j-\ell}}{(a)_i(a)_j}.
\]
On the other hand, the right-hand side of \eqref{eq:starformula} becomes
\[
\sum_{k=0}^{\min(i,j)}\frac{i^{\underline{k}}j^{\underline{k}}}
{k!(a)_k}y^{i+j-k}(1+y)^k.
\]
Its coefficient of $y^{i+j-\ell}$ is
\[
B_\ell=\sum_{k=\ell}^{\min(i,j)}
\frac{i^{\underline{k}}j^{\underline{k}}}{k!(a)_k}\binom{k}{\ell}.
\]
Put $I=i-\ell$, $J=j-\ell$, and $k=\ell+r$. Splitting the rising and
falling factorials at $\ell$ yields
\[
B_\ell=\frac{i^{\underline{\ell}}j^{\underline{\ell}}}
{\ell!(a)_\ell}\sum_{r=0}^{\min(I,J)}
\frac{I^{\underline{r}}J^{\underline{r}}}{r!(a+\ell)_r}
=\frac{i^{\underline{\ell}}j^{\underline{\ell}}}
{\ell!(a)_\ell}\cdot\frac{(a+j)_{i-\ell}}{(a+\ell)_{i-\ell}},
\]
where the last equality uses Lemma~\ref{lem:pochsum} with $b=a+\ell$.
Now
\[
\frac{i^{\underline{\ell}}j^{\underline{\ell}}}{\ell!}
=\binom{i}{\ell}\binom{j}{\ell}\ell!,\quad
(a)_\ell(a+\ell)_{i-\ell}=(a)_i,\quad
(a+j)_{i-\ell}=\frac{(a)_{i+j-\ell}}{(a)_j}.
\]
Thus $B_\ell=A_\ell$ for every $\ell$. This completes the proof.
\end{proof}

\section{Preservation of real-rootedness}\label{sec:rr}
We now prove Theorem~\ref{thm:rr-main}. Br\"and\'en and Ottergren
\cite{BO14} classified multiplier sequences for generalized Laguerre
bases, establishing a powerful criterion for real-rootedness
preservation in that setting. Remarkably, our weighted Euler numerators
map precisely onto their characterizing polynomial, which allows us to
invoke their classification directly and bypass the standard technical
estimates. (For completeness, a direct proof using the classical
P\'olya--Schur criterion and Br\"and\'en's bilinear theorem is given in
Appendix~\ref{app:rr}.)

\subsection{The change of variables and the location of zeros}

For a polynomial \(P\) of degree at most \(m\), define
\[
f(y)=(1+y)^mP\left(\frac{y}{1+y}\right).
\]
It is straightforward to check that \(f\) is a polynomial: writing
\(P(x)=\sum_{j=0}^mp_jx^j\) yields
\(f(y)=\sum_{j=0}^mp_jy^j(1+y)^{m-j}\). The inverse substitution is
\[
P(x)=(1-x)^mf\left(\frac{x}{1-x}\right).
\]
In particular, applying \eqref{eq:Lint} gives
\begin{equation}\label{eq:L-H}
\mathcal{L}_a(F)(y)=(1+y)^mH_{a,F}\left(\frac{y}{1+y}\right),\qquad
m=\deg F.
\end{equation}

\begin{lemma}\label{lem:zerochange}
Let \(P\in\mathbb{R}[x]\) be a nonzero polynomial with \(\deg P\le m\)
for some \(m\ge0\). If every zero of \(P\) is real and nonpositive,
then \(f(y)=(1+y)^mP(\frac{y}{1+y})\) has degree \(m\) and all its
zeros in \([-1,0]\). Conversely, if a nonzero polynomial \(f\) has
degree \(m\) and all its zeros in \([-1,0]\), then
\(P(x)=(1-x)^mf(\frac{x}{1-x})\) is a nonzero polynomial of degree at
most \(m\) with only real, nonpositive zeros. Under this
correspondence, we have \(P(0)=f(0)\); hence \(P(0)\neq0\) if and only
if \(f(0)\neq0\). Moreover, if \(f\) has positive leading coefficient,
then \(P\) has nonnegative coefficients.
\end{lemma}
\begin{proof}
Write \(P(x)=c\prod_{\nu=1}^d(x+r_\nu)\) with \(c\neq0\) and
\(r_\nu\ge0\), where \(d=\deg P\). Then
\[
f(y)=c(1+y)^{m-d}\prod_{\nu=1}^d((1+r_\nu)y+r_\nu).
\]
This has degree \(m\), with roots at \(-r_\nu/(1+r_\nu)\in(-1,0]\) for
each \(\nu\), and at \(-1\) with multiplicity \(m-d\).

Conversely, suppose \(f(y)=c\prod_{\nu=1}^m(y+\rho_\nu)\) with
\(c\neq0\) and \(0\le\rho_\nu\le1\). The inverse substitution yields
\[
P(x)=c\prod_{\nu=1}^m(\rho_\nu+(1-\rho_\nu)x).
\]
When \(\rho_\nu<1\), this contributes a nonpositive root
\(-\rho_\nu/(1-\rho_\nu)\); when \(\rho_\nu=1\), the factor is the
constant \(1\). Therefore \(P\) is nonzero, and its degree is \(m\)
minus the multiplicity of \(-1\) as a root of \(f\). The equality
\(P(0)=f(0)\) gives the claim about the origin. Finally, all factors
in the displayed expression for \(P\) have nonnegative coefficients, so
the same holds for \(P\) when \(c>0\).
\end{proof}

\subsection{Multiplier sequences for the Laguerre basis}

Fix \(a>0\) and put \(\alpha=a-1>-1\). Let
\[
L_n^{(\alpha)}(u)=\sum_{j=0}^n\frac{(-1)^j}{j!}
\binom{n+\alpha}{n-j}u^j,\qquad n\ge0,
\]
be the generalized Laguerre polynomials. Since their degree is exactly
\(n\), they form a basis of \(\mathbb{R}[u]\).

For a nonzero polynomial \(F\in\mathbb{R}[t]\), define the linear operator
\(T_F\) on \(\mathbb{R}[u]\) by
\[
T_F\bigl(L_n^{(\alpha)}(u)\bigr)=F(n)L_n^{(\alpha)}(u),\qquad n\ge0.
\]
(Here and below, preservation of real-rootedness includes mapping to
the zero polynomial.)

We use the following characterization due to Br\"and\'en and Ottergren
\cite[Theorem 1.1]{BO14}, stated here in a form adapted to our notation.
Given real sequences \(\{a_j\}_{j\ge0}\) and \(\{\lambda_n\}_{n\ge0}\)
related by
\[
\lambda_n = \sum_{j=0}^n a_j \binom{n}{j},
\]
assume that \(\{\lambda_n\}\) is not a \emph{trivial sequence} (i.e., it
does not consist of at most two consecutive nonzero terms). Define the
power series
\[
p(y) = \sum_{j=0}^\infty \binom{j+\alpha}{j} a_j y^j.
\]
Then the operator \(T(L_n^{(\alpha)}(u)) = \lambda_n L_n^{(\alpha)}(u)\)
preserves real-rootedness if and only if \(p(y)\) is a real-rooted
polynomial with all its zeros contained in \([-1,0]\).

We now apply this to the sequence \(\lambda_n = F(n)\). Writing
\(F(t)=\sum_{j=0}^m c_jt^{\underline{j}}\) and using
\(t^{\underline{j}}=j!\binom{t}{j}\), we obtain
\[
F(n)=\sum_{j=0}^m (j!c_j)\binom{n}{j}.
\]
Thus, the sequence corresponds to \(a_j = j!c_j\) for \(0\le j\le m\),
and \(a_j=0\) for \(j>m\). Since \(F\) is a nonzero polynomial, it has
at most \(m\) roots, so \(F(n)\neq 0\) for infinitely many \(n\). Hence
the sequence \((F(n))_{n\ge0}\) is not trivial. For this sequence, the
power series \(p(y)\) truncates to the polynomial
\begin{equation}\label{eq:BO-poly}
p(y)
=\sum_{j=0}^m\binom{j+\alpha}{j}(j!c_j)y^j
=\sum_{j=0}^m\binom{j+a-1}{j}(j!c_j)y^j
=\sum_{j=0}^mc_j(a)_jy^j.
\end{equation}

Comparing this with \eqref{eq:La}, we see that \(p(y)\) is precisely
\(\mathcal{L}_a(F)(y)\). Therefore, the theorem of
\cite[Theorem 1.1]{BO14} gives
\begin{equation}\label{eq:BO14-equiv}
T_F\text{ preserves real-rootedness}
\iff
\mathcal{L}_a(F)\text{ has all its zeros in }[-1,0].
\end{equation}

By Lemma~\ref{lem:zerochange} and \eqref{eq:L-H}, the right-hand side
is equivalent to saying that all zeros of \(H_{a,F}\) are real and nonpositive.
We conclude that
\begin{equation}\label{eq:T-H}
T_F\text{ preserves real-rootedness}
\iff
H_{a,F}\text{ has only real nonpositive zeros}.
\end{equation}

Now we are ready to prove the real-rootedness of the Euler numerators.
\begin{proof}[Proof of Theorem~\ref{thm:rr-main}]
Assume that \(H_{a,F}\) and \(H_{a,G}\) have only real nonpositive
zeros. By \eqref{eq:T-H}, both \(T_F\) and \(T_G\) preserve
real-rootedness; consequently, their composition \(T_F\circ T_G\) also
preserves real-rootedness. For every \(n\ge0\), we compute
\[
(T_F\circ T_G)\bigl(L_n^{(\alpha)}\bigr)
=T_F\bigl(G(n)L_n^{(\alpha)}\bigr)
=F(n)G(n)L_n^{(\alpha)}
=T_{FG}\bigl(L_n^{(\alpha)}\bigr).
\]
Since these two linear operators agree on the basis
\(\{L_n^{(\alpha)}\}\), they must be identical. Thus,
\(T_{FG}\) also preserves real-rootedness. By applying the equivalence
\eqref{eq:T-H} once more, we deduce that all zeros of \(H_{a,FG}\)
are real and nonpositive.
\end{proof}

\section{Preservation of \(\gamma\)-positivity}\label{sec:gamma}
We now focus on Theorem~\ref{thm:gamma-main}. We begin by fixing
terminology. A nonzero polynomial \(P\in\mathbb{R}[x]\) is
\emph{symmetric of degree} \(m\in\mathbb{Z}_{\ge0}\) if
\[
\deg P\le m\qquad\text{and}\qquad P(x)=x^mP(1/x).
\]
Equivalently, this means
\(p_j=p_{m-j}\) for every \(j\) when writing \(P(x)=\sum_{j=0}^mp_jx^j\). 
If \(J=\operatorname{supp}(P)\)
denotes the set of indices with nonzero coefficients, then such an \(m\)
is unique and equals \(\min J+\max J\). We call this the
\emph{symmetric degree} of \(P\) and denote it by
\[
\operatorname{sdeg}(P)=\min\{j:p_j\neq0\}+\max\{j:p_j\neq0\}.
\]

The following basic properties follow directly from the definition:
\begin{enumerate}
\item[(S1)] \(\operatorname{sdeg}(1)=0\),
\(\operatorname{sdeg}(1+x)=1\), \(\operatorname{sdeg}(x)=2\), and
\(\operatorname{sdeg}(x^n)=2n\) for \(n\ge0\).
\item[(S2)] \(\operatorname{sdeg}(PQ)=\operatorname{sdeg}(P)+
\operatorname{sdeg}(Q)\) for nonzero \(P,Q\in\mathbb{R}[x]\).
\item[(S3)] \(\operatorname{sdeg}(P)=\deg P\) if and only if
\(p_0\neq0\).
\end{enumerate}
By (S3), the symmetric degree can exceed the ordinary degree exactly
when the constant term vanishes. For instance, \(\deg(ax)=1\) but
\(\operatorname{sdeg}(ax)=2\) for any \(a\neq0\).

Recall from \eqref{eq:Gammam} that \(\Gamma_m\) is the cone generated by
the polynomials \(x^i(1+x)^{m-2i}\). 

The polynomials
\begin{equation}\label{eq:Bmi}
B_{m,i}(x)=x^i(1+x)^{m-2i},\qquad0\le i\le\lfloor m/2\rfloor,
\end{equation}
form a basis of the space of polynomials symmetric of degree \(m\).
Indeed, the space has dimension \(\lfloor m/2\rfloor+1\), and the
degrees \(\deg B_{m,i}=m-i\) are all distinct for \(0 \le i \le \lfloor m/2\rfloor\). Therefore, every such
\(P\) has a unique expansion
\begin{equation}\label{eq:gammaexp}
P(x)=\sum_{i=0}^{\lfloor m/2\rfloor}\gamma_ix^i(1+x)^{m-2i}.
\end{equation}
Consequently, \(P\in\Gamma_m\) if and only if \(\operatorname{sdeg}(P)=m\)
and all coefficients \(\gamma_i\) in this expansion are nonnegative.
Equivalently, \(\Gamma_m\) is precisely the cone generated by the
\(B_{m,i}\).

We will use the following elementary properties:
\begin{enumerate}
\item[(G1)] \(x^k(1+x)^\ell\in\Gamma_{2k+\ell}\) for all
\(k,\ell\ge0\). In particular, \(x\in\Gamma_2\) and
\(1+x\in\Gamma_1\).
\item[(G2)] \(\Gamma_m\cdot\Gamma_n\subseteq\Gamma_{m+n}\): if
\(P\in\Gamma_m\) and \(Q\in\Gamma_n\), then \(PQ\in\Gamma_{m+n}\).
Indeed, \(B_{m,i}B_{n,j}=B_{m+n,i+j}\), so multiplying the
\(\gamma\)-expansions \eqref{eq:gammaexp} gives a nonnegative
expansion in the basis of \(\Gamma_{m+n}\).
\item[(G3)] For \(\alpha>0\),
\((x+\alpha)(x+\alpha^{-1})\in\Gamma_2\). Indeed,
\[
(x+\alpha)(x+\alpha^{-1})=x^2+(\alpha+\alpha^{-1})x+1
=(1+x)^2+(\alpha+\alpha^{-1}-2)x,
\]
and the coefficient of \(x\) is nonnegative.
\end{enumerate}
Note that (G3) is the only step where the reciprocal pairing of roots
appears. Combining (G1)--(G3) with the factorization of real-rooted
polynomials yields the following classical result, which is due to Br\"and\'en
\cite[Lemma 4.1]{Branden2008} and Gal \cite[Remark 3.1.1]{Gal2005}; see also
\cite[Section~1]{Athanasiadis}.

\begin{theorem}[Br\"and\'en--Gal]\label{thm:BrandenGal}
Let \(P\in\mathbb{R}_{\ge0}[x]\) be symmetric of degree \(m\). If all
zeros of \(P\) are real and nonpositive, then \(P\in\Gamma_m\).
\end{theorem}

Next, we verify that symmetry is preserved under our operation.
\begin{proposition}\label{prop:sym}
Let \(a>0\), and let \(F,G\in\mathbb{R}[t]\) be nonzero polynomials of
degrees \(m,n\) respectively. If \(H_{a,F}\) and \(H_{a,G}\) are
symmetric of degrees \(m\) and \(n\), then \(H_{a,FG}\) is symmetric
of degree \(m+n\).
\end{proposition}

\begin{proof}
Put \(P=H_{a,F}\) and \(f=\mathcal{L}_a(F)\). By \eqref{eq:Hint}, we
have
\begin{equation}\label{eq:finv}
f(y)=(1+y)^mP\left(\frac{y}{1+y}\right).
\end{equation}
We claim that \(P\) is symmetric of degree \(m\) if and only if
\begin{equation}\label{eq:symequiv}
f(-1-y)=(-1)^mf(y).
\end{equation}
For the forward implication, suppose \(P(x)=x^mP(1/x)\). Substituting \(y=-1-x\) into
\eqref{eq:finv} yields
\begin{align*}
f(-1-x) &= (-x)^mP\left(\frac{-1-x}{-x}\right) \\
&= (-x)^mP\left(\frac{1+x}{x}\right) \\
&= (-x)^m\left(\frac{1+x}{x}\right)^mP\left(\frac{x}{1+x}\right)
= (-1)^mf(x).
\end{align*}
The backward implication is obtained similarly using the substitution
\(y=\frac{x}{1-x}\).

Analogously, for \(g=\mathcal{L}_a(G)\), we have
\(g(-1-y)=(-1)^ng(y)\). Differentiating \(k\) times gives
\[
f^{(k)}(-1-y)=(-1)^{m+k}f^{(k)}(y),\qquad
g^{(k)}(-1-y)=(-1)^{n+k}g^{(k)}(y).
\]
Since \((-1-y)(-y)=y(y+1)\), each summand in \eqref{eq:starformula}
acquires the same factor \((-1)^{m+n}\) under \(y\mapsto-1-y\). Hence
\[
(f\star_ag)(-1-y)=(-1)^{m+n}(f\star_ag)(y).
\]
Because \(f\star_ag=\mathcal{L}_a(FG)\), we apply the equivalence
\eqref{eq:symequiv} once again to conclude that \(H_{a,FG}\) is
symmetric of degree \(m+n\).
\end{proof}

We now prove Theorem~\ref{thm:gamma-main}. Note that the two
gamma-positive inputs are not required to be real-rooted.

\begin{proof}[Proof of Theorem~\ref{thm:gamma-main}]
Write
\[
H_{a,F}(x)=\sum_i\gamma_iB_{m,i}(x),\qquad
H_{a,G}(x)=\sum_j\eta_jB_{n,j}(x),
\]
with all \(\gamma_i,\eta_j\ge0\).

For each \(0\le i\le\lfloor m/2\rfloor\), define
\begin{equation}\label{eq:fi}
f_i(y)=(1+y)^mB_{m,i}\left(\frac{y}{1+y}\right)
=[y(1+y)]^i(1+2y)^{m-2i},\qquad
F_i=\mathcal{L}_a^{-1}(f_i).
\end{equation}
Although \(B_{m,i}\) has degree \(m-i\), the polynomial \(f_i\) has
degree exactly \(m\), with leading coefficient \(2^{m-2i}\). Since
\(\mathcal{L}_a\) preserves degrees,
\[
\deg F_i=m,\qquad\lc(F_i)=\frac{2^{m-2i}}{(a)_m}>0,\qquad
H_{a,F_i}=B_{m,i}.
\]
The polynomials \(G_j\) are defined analogously, and we have
\[
\deg G_j=n,\qquad\lc(G_j)=\frac{2^{n-2j}}{(a)_n}>0,\qquad
H_{a,G_j}=B_{n,j}.
\]

Fix \(i,j\) and set \(U_{ij}=H_{a,F_iG_j}\). We show that
\(U_{ij}\in\Gamma_{m+n}\) by verifying the conditions of the
Br\"and\'en--Gal theorem.

First, since \(B_{m,i}\) and \(B_{n,j}\) have no zeros outside
\(\{0,-1\}\), Theorem~\ref{thm:rr-main} implies that all zeros of
\(U_{ij}\) are real and nonpositive.

Second, by \eqref{eq:Hvalues} of Proposition 2.1, we have
\[
U_{ij}(1)=H_{a,F_iG_j}(1)=(a)_{m+n}\lc(F_i)\lc(G_j)>0.
\]
Since \(U_{ij}\) has only real nonpositive zeros, being positive at \(1\)
ensures that all its coefficients are nonnegative.

Third, by Proposition~\ref{prop:sym}, \(U_{ij}\) is symmetric of degree
\(m+n\).

Theorem~\ref{thm:BrandenGal} therefore applies to yield
\(U_{ij}\in\Gamma_{m+n}\).

Finally, by the linearity of \(\mathcal{L}_a^{-1}\) and \eqref{eq:fi},
we have
\[
F=\sum_i\gamma_iF_i,\qquad G=\sum_j\eta_jG_j.
\]
All products \(F_iG_j\) have degree \(m+n\), so their weighted Euler
series share the same denominator \((1-x)^{a+m+n}\). Consequently,
\[
H_{a,FG}(x)=\sum_{i,j}\gamma_i\eta_jU_{ij}(x).
\]
The theorem then follows, since the cone \(\Gamma_{m+n}\) is closed under nonnegative linear
combinations.
\end{proof}

\begin{remark}\label{rem:subtle-degree}
The requirement that \(\deg F = m\) and \(\deg G = n\) in
Theorem~\ref{thm:gamma-main} is essential. The map \(F \mapsto H_{a,F}\)
is not globally linear: the denominator of its generating function
depends on \(\deg F\). For example, with \(a=1\), \(H_{1,t}=x\) and
\(H_{1,1}=1\), but
\[
H_{1,t+1}=1 \neq x+1=H_{1,t}+H_{1,1}.
\]
By imposing \(\deg F = m\) and \(H_{a,F}\in\Gamma_m\), the ordinary
degree and the symmetric degree are forced to coincide. This unifies
the denominators, ensuring that the linear combinations used in the
proof (such as \(F=\sum_i\gamma_iF_i\)) are well-defined.

Without this degree condition, the conclusion can fail. For instance,
take \(a=1\) and \(F(t)=G(t)=(t+1)^2\). Then 
\(H_{1,F}=H_{1,G}=1+x\in\Gamma_1\). Direct computation yields
\(
H_{1,FG}=1+11x+11x^2+x^3\in\Gamma_3,
\)
which does not belong to the predicted  
\(\Gamma_{1+1}=\Gamma_2\). This shows that the degree condition is indeed necessary.
\end{remark}

\section{A common simplex factor}\label{sec:geom}

\subsection{Ehrhart series and the Pochhammer weight}
Let \(P\) be a \(d\)-dimensional lattice polytope. We work in its affine
span and, after translating by a lattice vertex, regard \(P\) as a
full-dimensional polytope in \(\mathbb{R}^d\) with the induced lattice;
the lattice on a Cartesian product is the product lattice (see
\cite{BeckRobins} for standard background). The Ehrhart function of \(P\) is
\[
L_P(k)=|kP\cap\mathbb{Z}^d|,\qquad k\in\mathbb{Z}_{\ge0}.
\]
Ehrhart's theorem states that \(L_P(k)\) is a polynomial in \(k\) of degree
\(d\), with \(L_P(0)=1\) \cite{Ehrhart}, and its generating function takes
the form
\begin{equation}\label{eq:ehrhart}
\sum_{k\ge0}L_P(k)x^k=\frac{h_P^*(x)}{(1-x)^{d+1}}.
\end{equation}
The numerator \(h_P^*(x)\) is called the \(h^*\)-polynomial of \(P\). It has
degree at most \(d\), nonnegative integer coefficients, and constant term
\(1\); see \cite{Stanley,BeckRobins}. In particular, if \(h_P^*\) is
real-rooted, its zeros are strictly negative. (As usual, nonzero
constants are considered real-rooted.)

For two lattice polytopes \(P\) and \(Q\), the Cartesian product satisfies
\begin{equation}\label{eq:ehrhartprod}
L_{P\times Q}(k)=L_P(k)L_Q(k).
\end{equation} 
The preservation of real-rootedness
under this product is a classical result of Wagner
\cite[Theorem 0.2]{Wagner}; see Liu \cite[Section 2]{Liu} for a broader
survey of operations on lattice polytopes.

Fix an integer \(r\ge1\) and let
\[
\Delta_{r-1}=\operatorname{conv}(0,e_1,\dots,e_{r-1})\subseteq
\mathbb{R}^{r-1}
\]
be the standard simplex, with \(\Delta_0\) being a point. The lattice
points of \(k\Delta_{r-1}\) correspond to nonnegative integer tuples
\((u_1,\dots,u_{r-1})\) satisfying \(u_1+\dots+u_{r-1}\le k\). By
setting \(u_r=k-u_1-\dots-u_{r-1}\), these are in bijection with the
nonnegative integer solutions of \(u_1+\dots+u_r=k\). Consequently,
\begin{equation}\label{eq:simplex}
L_{\Delta_{r-1}}(k)=\binom{k+r-1}{r-1}=\frac{(r)_k}{k!}.
\end{equation}
Equivalently, the polynomial \(W_r(t)=\binom{t+r-1}{r-1}\) is the
Ehrhart polynomial of \(\Delta_{r-1}\) (with the empty product being
\(1\) for \(r=1\)). The Pochhammer weight \((a)_k/k!\) therefore has a
direct geometric interpretation at every positive integer parameter.

\begin{proposition}\label{prop:H-ehrhart}
Let \(P\) be a lattice polytope of dimension \(d\), and let \(r\ge1\) be
an integer. Then
\begin{equation}\label{eq:H-ehrhart}
H_{r,L_P}(x)=h^*_{P\times\Delta_{r-1}}(x).
\end{equation}
In particular, \(\deg h^*_{P\times\Delta_{r-1}}\le d\).
\end{proposition}

\begin{proof}
By \eqref{eq:ehrhartprod} and \eqref{eq:simplex}, we have
\[
\sum_{k\ge0}\frac{(r)_k}{k!}L_P(k)x^k
=\sum_{k\ge0}L_{P\times\Delta_{r-1}}(k)x^k.
\]
The left-hand side is the weighted Euler series associated to \(F=L_P\).
Since \(\deg L_P=d\), its denominator is \((1-x)^{r+d}\). On the other
hand, the polytope \(P\times\Delta_{r-1}\) has dimension \(d+r-1\), so
its Ehrhart series shares the same denominator \((1-x)^{r+d}\).
Comparing the numerators yields \eqref{eq:H-ehrhart}. The degree bound
then follows from the degree bound for \(H_{r,L_P}\).
\end{proof}
\subsection{Products with one common simplex}

We now prove Theorems~\ref{thm:geom-rr} and~\ref{thm:geom-gamma}. Using
Proposition~\ref{prop:H-ehrhart} to translate into the algebraic
setting, both proofs are immediate applications of our main results
(Theorems~\ref{thm:rr-main} and~\ref{thm:gamma-main}). Afterward, we
explain a geometric feature of the output.

\begin{proof}[Proof of Theorems~\ref{thm:geom-rr} and~\ref{thm:geom-gamma}]
Put \(F=L_P\) and \(G=L_Q\). Then
\(\deg F=d\), \(\deg G=e\), and \(L_{P\times Q}=FG\).
By Proposition~\ref{prop:H-ehrhart}, we have
\[
H_{r,F}=h^*_{P\times\Delta_{r-1}},\qquad
H_{r,G}=h^*_{Q\times\Delta_{r-1}},\qquad
H_{r,FG}=h^*_{P\times Q\times\Delta_{r-1}}.
\]

Under the hypotheses of Theorem~\ref{thm:geom-rr}, the first two
numerators are real-rooted. Since Ehrhart numerators have nonnegative
coefficients and constant term \(1\), all their zeros are strictly
negative. Applying Theorem~\ref{thm:rr-main} with \(a=r\) shows that
\(H_{r,FG}\) has only real nonpositive zeros, proving
Theorem~\ref{thm:geom-rr}.

Under the hypotheses of Theorem~\ref{thm:geom-gamma}, we instead have
\(H_{r,F}\in\Gamma_d\) and \(H_{r,G}\in\Gamma_e\). Together with
\(\deg F=d\) and \(\deg G=e\), these are precisely the hypotheses of
Theorem~\ref{thm:gamma-main} with \(a=r\). Hence
\[
H_{r,FG}\in\Gamma_{d+e},
\]
which proves Theorem~\ref{thm:geom-gamma}.
\end{proof}

\begin{corollary}\label{cor:closure}
For each positive integer \(r\), the class
\[
\mathcal{C}_r=\{P:h^*_{P\times\Delta_{r-1}}\text{ is real-rooted}\}
\]
is closed under Cartesian products. More generally, if
\(P_1,\dots,P_s\in\mathcal{C}_r\), then
\(h^*_{P_1\times\dots\times P_s\times\Delta_{r-1}}\) is real-rooted.
Likewise, if \(P_1,\dots,P_s\) are lattice polytopes with
\(h^*_{P_i\times\Delta_{r-1}}\in\Gamma_{d_i}\) for every \(i\), where
\(d_i=\dim P_i\), then the same output belongs to
\(\Gamma_{d_1+\dots+d_s}\).
\end{corollary}

\begin{proof}
The assertion for two factors is the combination of
Theorem~\ref{thm:geom-rr} and Theorem~\ref{thm:geom-gamma}. Apply them
successively to \(P_1\times\dots\times P_j\) and \(P_{j+1}\),
keeping the same parameter \(r\). At the \(j\)th step the symmetry degree in
the second assertion is \(d_1+\dots+d_j\).
\end{proof}

\subsection{A common Gorenstein index}
We next discuss a geometric condition that supplies the symmetry degrees
required in Theorem~\ref{thm:geom-gamma}. A full-dimensional lattice
polytope \(A\subseteq\mathbb{R}^D\) is \emph{reflexive} if \(0\) lies in
its interior and its polar
\[
A^\vee=\{y\in\mathbb{R}^D:\langle y,x\rangle\le1
\text{ for every }x\in A\}
\]
is a lattice polytope. Equivalently, every facet inequality can be
written as \(\langle n,x\rangle\le1\) for some primitive integer vector
\(n\).

The \emph{codegree} of a lattice polytope \(P\) is
\[
\operatorname{codeg}P=\min\{k\ge1:\operatorname{relint}(kP)\cap
\mathbb{Z}^d\neq\emptyset\}.
\]
We say \(P\) is \emph{Gorenstein of index \(r\)} if
\(r=\operatorname{codeg}P\), the polytope \(rP\) has a unique interior
lattice point \(u\), and \(rP-u\) is reflexive. These definitions are
taken in the affine lattice of \(P\). A point is Gorenstein of index
\(1\). We refer to \cite{BrunsRomer,BeckRobins} for the standard
connections with Ehrhart series.

We will need the following well-known properties of Gorenstein polytopes;
see \cite{BeckRobins} and \cite[Theorem 4.3]{LTX2026}.

\begin{lemma}\label{lem:gorenstein}
Let \(A\) be a \(D\)-dimensional Gorenstein lattice polytope of index
\(r\). Then
\begin{equation}\label{eq:gorenstein}
h_A^*(x)=x^{D+1-r}h_A^*(1/x),\qquad\deg h_A^*=D+1-r.
\end{equation}
Moreover, the Cartesian product of two Gorenstein lattice polytopes of
the same index \(r\) is Gorenstein of index \(r\).
\end{lemma}

For \(r\ge 2\), the polytope \(r\Delta_{r-1}\) has a unique interior
lattice point, namely the vector \(\mathbf{1}=(1,\dots,1)\in\mathbb{Z}^{r-1}\).
Indeed, the interior of \(k\Delta_{r-1}\) consists of points with
strictly positive coordinates and coordinate sum less than \(k\). Since
\(\mathbf{1}\) has coordinate sum \(r-1 < r\), it lies in the interior
of \(r\Delta_{r-1}\); moreover, no such point exists for \(k<r\).
Translating by \(-\mathbf{1}\), the facet inequalities become
\[
-x_i\le1\quad(1\le i\le r-1),\qquad
x_1+\dots+x_{r-1}\le1.
\]
Since the right-hand sides are all \(1\), the polytope is reflexive.
Hence \(\Delta_{r-1}\) is Gorenstein of index \(r\). (For \(r=1\), this
holds trivially by our convention for a point.)

\begin{corollary}\label{cor:gorenstein}
Let \(P\) and \(Q\) be Gorenstein lattice polytopes of the same index \(r\).
If \(h^*_{P\times\Delta_{r-1}}\) and \(h^*_{Q\times\Delta_{r-1}}\) are
gamma-positive, then \(h^*_{P\times Q\times\Delta_{r-1}}\) is
gamma-positive.
\end{corollary}

\begin{proof}
Put \(d=\dim P\) and \(e=\dim Q\). Lemma~\ref{lem:gorenstein} shows that
\(P\times\Delta_{r-1}\) is Gorenstein of index \(r\). Its dimension is
\(d+r-1\), so its \(h^*\)-polynomial has degree and symmetry degree
\((d+r-1)+1-r=d\). Similarly, the symmetry degree for
\(Q\times\Delta_{r-1}\) is \(e\). The two gamma-positive inputs therefore
belong to \(\Gamma_d\) and \(\Gamma_e\), respectively.
Theorem~\ref{thm:geom-gamma} gives the conclusion, with symmetry degree
\(d+e\).
\end{proof}

\subsection{The Reeve tetrahedron}
Consider the Reeve tetrahedron \cite{Reeve}
\[
R=\operatorname{conv}\big(0,(1,0,0),(0,1,0),(1,1,2)\big).
\]
A direct computation of its Ehrhart series gives
\begin{equation}\label{eq:hR}
L_R(k) = \frac{(k+1)\left((k+1)^2+2\right)}{3}, \qquad
h_R^*(x)=1+x^2.
\end{equation}
This numerator is not real-rooted. We give two applications of our
theorems to this non-real-rooted example.

\paragraph{First application: adjoining a segment.}
Let \(\Delta_1\) be the unit segment. Since its Ehrhart polynomial is
\(k+1\), multiplying the \(k\)th series coefficient by \(k+1\) acts as
\(1+x\frac{d}{dx}\). We obtain
\begin{align*}
\sum_{k\ge0}L_{R\times\Delta_1}(k)x^k
&=\left(1+x\frac{d}{dx}\right)\frac{1+x^2}{(1-x)^4}\\
&=\frac{(1-x)(1+3x^2)+4x(1+x^2)}{(1-x)^5}
=\frac{(1+x)^3}{(1-x)^5}.
\end{align*}
Thus \(h^*_{R\times\Delta_1}(x)=(1+x)^3\). In particular, the
non-real-rooted numerator \(1+x^2\) is repaired by adjoining a single
segment: \(R\notin\mathcal{C}_1\) but \(R\in\mathcal{C}_2\). Moreover,
since the output belongs to \(\Gamma_3\),
Corollary~\ref{cor:closure} implies that for every \(s\ge1\), the
numerator \(h^*_{R^{\times s}\times\Delta_1}(x)\) is real-rooted and
belongs to \(\Gamma_{3s}\).

\paragraph{Second application: Cartesian powers.}
Next we consider the powers \(R^{\times s}\) without adjoining any
simplex factor. For \(s\ge1\), define
\begin{equation}\label{eq:Bs-def}
B_s(x) = h^*_{R^{\times s}}(x), \qquad
\sum_{k\ge0}L_R(k)^s x^k = \frac{B_s(x)}{(1-x)^{3s+1}}.
\end{equation}
We will show that \(B_s\) has only real negative roots if and only if
\(s\ge3\).

Since \(B_s(0)=1\), zero is never a root. For \(s=1\), \(B_1(x)=1+x^2\)
has no real roots. For \(s=2\), a direct computation gives
\[
B_2(x)=(1+x)(x^4+8x^3+22x^2+8x+1),
\]
and the quartic factor equals \((x^2+4x+1)^2+4x^2>0\) for all real \(x\).
Hence \(B_2\) has only one real root. Thus the threshold \(s\ge3\) is
sharp.

We now prove the assertion for \(s\ge3\). By Theorem~\ref{thm:rr-main}
with \(a=1\), the product of two real-rooted polynomials with real
nonpositive zeros has real nonpositive zeros. Since
\(L_{R^{\times p}}(k)=L_R(k)^p\), this implies
\begin{equation}\label{eq:Bs-product}
B_p,\,B_q \text{ real-rooted} \implies B_{p+q} \text{ real-rooted}.
\end{equation}
It therefore suffices to verify the base cases \(s=3,4,5\), as the
identities \(6=3+3\), \(7=3+4\), \(8=3+5\), and the induction step
\(s=3+(s-3)\) for \(s\ge6\) complete the proof.

To check these base cases, we make the substitution \(u=-x-1/x\). For
\(u>2\), the roots of \(x^2+ux+1=0\) are real and negative. We express
\(B_s\) in terms of a lower-degree monic polynomial \(q_s\). Put
\(d=3s-1\), \(\ell=\lfloor d/2\rfloor\), and \(\varepsilon=d-2\ell\).
Then
\[
B_s(x)=(-1)^\ell x^\ell(1+x)^\varepsilon q_s(-x-1/x).
\]
The extra factor \(1+x\) occurs when \(s\) is even. The coefficients
of \(B_s\) are given by \(L_R\) via the generating function.

For \(s=3\), \(B_3\) has degree \(8\):
\begin{align*}
B_3(x)={}&1+54x+736x^2+3274x^3+5310x^4+3274x^5+736x^6+54x^7+x^8.
\end{align*}
It can be written as \(B_3(x)=x^4 q_3(-x-1/x)\), with
\[
q_3(u)=u^4-54u^3+732u^2-3112u+3840.
\]
Evaluating \(q_3\) at \(u=2,3,5,12,36\) gives
\(128,-285,455,-672,672\). The signs alternate, so the intermediate
value theorem ensures a root in each interval
\((2,3), (3,5), (5,12), (12,36)\). All evaluations use exact rational
arithmetic. Since \(\deg q_3=4\), all its roots are real and lie in
\((2,\infty)\). Each root \(u_i\) yields two negative roots of \(B_3\)
via \(x^2+u_ix+1=0\). Thus all \(8\) roots of \(B_3\) are real and
negative.

For \(s=4,5\), the verification is analogous. The corresponding
polynomials and their alternating evaluations are:
\begin{align*}
B_4(x)&=-(1+x)x^5 q_4(-x-1/x), \\
q_4(u)&=u^5-242u^4+11144u^3-149008u^2+673616u-920480, \\
&\quad\text{evaluated at } u=2,3,5,12,38,187 \text{ yields } -, +, -, +, -, +.
\end{align*}
\begin{align*}
B_5(x)&=-x^7 q_5(-x-1/x), \\
q_5(u)&=u^7-1008u^6+144780u^5-5502080u^4+75156720u^3\\
&\quad -428237568u^2+1041166144u-892108800, \\
&\quad\text{evaluated at } u=2,\tfrac52,3,5,11,30,115,845 \text{ yields } -, +, -, +, -, +, -, +.
\end{align*}
Since \(\deg q_4=5\) and \(\deg q_5=7\), the alternating signs certify
that all their roots are real and lie in \((2,\infty)\). Translating
back via \(u=-x-1/x\) (and accounting for the extra factor \(1+x\) when
\(s=4\)) shows that all roots of \(B_4\) and \(B_5\) are real and
negative. This completes the proof of the base cases.

By \eqref{eq:Bs-product}, \(B_s\) has only real negative roots for all
\(s\ge3\).

\subsection*{Acknowledgements}
Guoce Xin was partially supported by the National Natural Science
Foundation of China (No. 12571355).

\appendix
\section{An alternative proof via the P\'olya--Schur criterion}\label{app:rr}

This appendix provides an alternative proof of
Theorem~\ref{thm:rr-main} that avoids the Br\"and\'en--Ottergren
characterization. Instead, it relies solely on the classical
P\'olya--Schur criterion and Br\"and\'en's bilinear preservation
theorem. The argument follows the structure of Section~\ref{sec:rr} but
substitutes the algebraic core.

\subsection{A Pochhammer multiplier sequence}

\begin{definition}\label{def:mult}
A sequence $(\lambda_k)_{k\ge0}$ of real numbers is a \emph{multiplier
sequence} if the coefficient operator
\[
T_\lambda\left(\sum_{k=0}^Nc_kz^k\right)=\sum_{k=0}^N\lambda_kc_kz^k
\]
maps every real-rooted polynomial to a real-rooted polynomial or to the
zero polynomial.
\end{definition}

We use the following finite form of the P\'olya--Schur criterion
\cite{PS}; see also \cite[Theorem 3.14]{B}.

\begin{theorem}[P\'olya--Schur criterion]\label{thm:PS}
Let $(\lambda_k)_{k\ge0}$ be a sequence of nonnegative real numbers. It
is a multiplier sequence if and only if, for every integer $N\ge0$,
the polynomial
\[
J_N(z)=\sum_{k=0}^N\binom{N}{k}\lambda_kz^k
\]
has only real, nonpositive zeros or is identically zero.
\end{theorem}

\begin{lemma}\label{lem:pochmult}
For every $a>0$, the sequence $(1/(a)_k)_{k\ge0}$ is a positive
multiplier sequence.
\end{lemma}

\begin{proof}
We verify Theorem~\ref{thm:PS} using generalized Laguerre polynomials.
For $\alpha>-1$ and $N\ge0$, define
\begin{equation}\label{eq:Laguerre}
L_N^{(\alpha)}(u)=\sum_{k=0}^N\frac{(-1)^k}{k!}
\binom{N+\alpha}{N-k}u^k.
\end{equation}
Here $\binom{v}{r}=v^{\underline{r}}/r!$ for $r\ge0$. Formula
\eqref{eq:Laguerre} is the standard finite expansion of the generalized
Laguerre polynomial \cite[Eq.~18.5.12]{NIST}. Its leading coefficient is
$(-1)^N/N!$, so its degree is exactly $N$.

To locate the zeros, we use the Rodrigues identity
\begin{equation}\label{eq:Rodrigues}
L_N^{(\alpha)}(u)=\frac{e^uu^{-\alpha}}{N!}\frac{d^N}{du^N}
\big(e^{-u}u^{N+\alpha}\big),\qquad u>0.
\end{equation}
Indeed, differentiating the exponential factor \(k\) times gives
\[
\frac{d^N}{du^N}\big(e^{-u}u^{N+\alpha}\big)
=e^{-u}\sum_{k=0}^N\binom{N}{k}(-1)^k
(N+\alpha)^{\underline{N-k}}u^{\alpha+k},
\]
and multiplying by \(e^uu^{-\alpha}/N!\) recovers \eqref{eq:Laguerre}.

Let $p\in\mathbb{R}[u]$ have degree less than $N$. Integrating by parts
$N$ times in \eqref{eq:Rodrigues} gives
\begin{equation}\label{eq:orth}
\int_0^\infty L_N^{(\alpha)}(u)p(u)u^\alpha e^{-u}\,du
=\frac{(-1)^N}{N!}\int_0^\infty p^{(N)}(u)e^{-u}u^{N+\alpha}\,du=0.
\end{equation}
To justify this, put $q(u)=e^{-u}u^{N+\alpha}$. For \(0\le r\le N-1\),
\(q^{(r)}(u)=O(u^{N+\alpha-r})\) as \(u\to0^+\), and
\(N+\alpha-r\ge\alpha+1>0\). Thus all boundary terms at \(0\) vanish.
At infinity, $q$ and its derivatives decay exponentially, so the other
boundary terms also vanish. All integrals converge since the least
power at \(0\) is \(u^\alpha\), with $\alpha>-1$.

Assume $N\ge1$, and let $u_1,\dots,u_s$ be the distinct zeros in
\((0,\infty)\) at which $L_N^{(\alpha)}$ changes sign. If $s<N$, take
\(p(u)=\prod_{\nu=1}^s(u-u_\nu)\). The polynomial $L_N^{(\alpha)}p$
then has constant sign on \((0,\infty)\) apart from its zeros, and it is
not identically zero. Since $u^\alpha e^{-u}>0$ on this interval, its
weighted integral cannot vanish. This contradicts \eqref{eq:orth},
because $\deg p=s<N$. Consequently $s\ge N$. As $L_N^{(\alpha)}$ has
degree $N$, it has exactly $N$ distinct positive zeros, each of
multiplicity one.

Now set $\alpha=a-1$. For \(0\le k\le N\),
\[
\binom{N+a-1}{N-k}=\frac{(a)_N}{(a)_k(N-k)!}.
\]
It follows from \eqref{eq:Laguerre} that
\begin{equation}\label{eq:Laguerre-special}
\sum_{k=0}^N\binom{N}{k}\frac{z^k}{(a)_k}
=\frac{N!}{(a)_N}L_N^{(a-1)}(-z).
\end{equation}
For $N\ge1$ all zeros of this polynomial are strictly negative; for
$N=0$ it equals \(1\). Theorem~\ref{thm:PS} proves the multiplier
property, and positivity follows from \(a>0\).
\end{proof}

\subsection{Interval preservation for the induced product}

We use the interval form of Br\"and\'en's bilinear preservation theorem
\cite[Theorem 3.11]{B}.

\begin{theorem}[Br\"and\'en]\label{thm:Branden}
Let $(\lambda_k)_{k\ge0}$ be a nonnegative multiplier sequence and let
$\alpha<\beta$ be real numbers. If $f\in\mathbb{R}[y]$ is real-rooted
and all zeros of $g\in\mathbb{R}[y]$ lie in $[\alpha,\beta]$, then
\[
\sum_{k\ge0}\frac{\lambda_k}{k!}(y-\alpha)^k(y-\beta)^k
f^{(k)}(y)g^{(k)}(y)
\]
is real-rooted or identically zero. (The sum is finite, and the closed
interval permits zeros at its endpoints.)
\end{theorem}

\begin{proposition}\label{prop:interval}
Fix $a>0$. Let $f,g\in\mathbb{R}[y]$ be nonzero polynomials with
positive leading coefficients, of degrees $m,n$, respectively. If all
their zeros belong to $[-1,0]$, then \(f\star_ag\) has degree $m+n$,
positive leading coefficient, and all its zeros in $[-1,0]$. Moreover,
\begin{equation}\label{eq:endpoint}
(f\star_ag)(0)=f(0)g(0),\qquad (f\star_ag)(-1)=f(-1)g(-1).
\end{equation}
\end{proposition}

\begin{proof}
By \eqref{eq:starformula},
\begin{equation}\label{eq:hstar}
h(y):=(f\star_ag)(y)=\sum_{k=0}^{\min(m,n)}
\frac{[y(y+1)]^k}{k!(a)_k}f^{(k)}(y)g^{(k)}(y).
\end{equation}
Lemma~\ref{lem:pochmult} allows us to apply
Theorem~\ref{thm:Branden} with
\[
\lambda_k=\frac{1}{(a)_k},\qquad\alpha=-1,\qquad\beta=0.
\]
Here $f$ is real-rooted, $g$ has all its zeros in $[-1,0]$, and
$(y-\alpha)^k(y-\beta)^k=[y(y+1)]^k$. The theorem therefore shows that
$h$ is real-rooted or identically zero.

Every summand in \eqref{eq:hstar} has degree at most $m+n$. Its
coefficient of $y^{m+n}$ is nonnegative, and their sum is
\begin{equation}\label{eq:lcstar}
\lc(f)\lc(g)\sum_{k=0}^{\min(m,n)}\frac{m^{\underline{k}}n^{\underline{k}}}
{k!(a)_k}>0.
\end{equation}
Thus $h$ is nonzero, has degree $m+n$, and has positive leading
coefficient.

It remains to locate its real zeros. Write
\(f(y)=c\prod_{\nu=1}^m(y+\rho_\nu)\), where \(c>0\) and \(0\le\rho_\nu\le1\).
For \(0\le k\le m\) and \(y>0\), each term in the product-rule expansion
of \(f^{(k)}(y)\) is positive. For \(y<-1\), each such term has sign
\((-1)^{m-k}\). The same statements hold for \(g\).

For \(y>0\), all summands in \eqref{eq:hstar} are positive. For \(y<-1\),
their signs are all \((-1)^{m-k+n-k}=(-1)^{m+n}\), since \(y(y+1)>0\).
In either interval the summand with \(k=0\) is nonzero, so \(h\) has no
zeros there. Together with real-rootedness, this proves that all its
zeros lie in \([-1,0]\). The same argument applies when \(m=0\) or
\(n=0\). Finally, evaluating \eqref{eq:hstar} at \(0\) and \(-1\) leaves
only its \(k=0\) summand and proves \eqref{eq:endpoint}.
\end{proof}

\subsection{Real-rootedness of the Euler numerators}

\begin{proof}[Alternative proof of Theorem~\ref{thm:rr-main}]
Set $m=\deg F$ and $n=\deg G$. Multiplying \(F\) and \(G\) by nonzero
real constants rescales their numerators and the numerator of \(FG\). We
may therefore assume that $P=H_{a,F}$ and $Q=H_{a,G}$ have positive
leading coefficients. Since all their zeros are nonpositive, their
coefficients are nonnegative, and $P(1),Q(1)>0$.

Set $f=\mathcal{L}_a(F)$ and $g=\mathcal{L}_a(G)$. By \eqref{eq:L-H}
and Lemma~\ref{lem:zerochange}, these polynomials have degrees $m,n$
and all their zeros in $[-1,0]$. Their leading coefficients are $P(1)$
and $Q(1)$, hence positive. The definition of the induced product gives
\[
\mathcal{L}_a(FG)=f\star_ag.
\]
Proposition~\ref{prop:interval} shows that this polynomial has degree
$m+n$, positive leading coefficient, and all its zeros in $[-1,0]$.
Since $\deg(FG)=m+n$, applying \eqref{eq:Hint} to \(FG\) and then
Lemma~\ref{lem:zerochange} proves that $H_{a,FG}$ has only real,
nonpositive zeros.

If both input numerators have only strictly negative zeros, then
$f(0)=P(0)\neq0$ and $g(0)=Q(0)\neq0$. By \eqref{eq:endpoint},
$(f\star_ag)(0)\neq0$, so the output numerator has no zero at the
origin. Its zeros are consequently strictly negative.

Finally, if the original input numerators already have nonnegative
coefficients, the initial rescaling is unnecessary. The positive
leading coefficient of \(f\star_ag\) and
Lemma~\ref{lem:zerochange} then show that \(H_{a,FG}\) has
nonnegative coefficients.
\end{proof}

\end{document}